\documentclass[10pt]{amsart}
\usepackage{amscd}
\usepackage{amssymb}

\newcommand{\md}{minimal diffeomorphism}
\newcommand{\Det}{{\mathrm{Det}}}

\theoremstyle{definition}
\newtheorem{thm}{Theorem}[section]
\newtheorem{lem}[thm]{Lemma}
\newtheorem{prp}[thm]{Proposition}
\newtheorem{dfn}[thm]{Definition}

\newtheorem{exa}[thm]{Example}

\newtheorem{qst}[thm]{Question}

\newcommand{\beq}{\begin{equation}}
\newcommand{\eeq}{\end{equation}}
\newcommand{\beqr}{\begin{eqnarray*}}
\newcommand{\eeqr}{\end{eqnarray*}}
\newcommand{\limi}[1]{\lim_{{#1} \to \infty}}
\newcommand{\bit}{\begin{itemize}}
\newcommand{\eit}{\end{itemize}}
\newcommand{\ts}[1]{{\textstyle{#1}}}

\newcommand{\af}{\alpha}
\newcommand{\bt}{\beta}
\newcommand{\gm}{\gamma}
\newcommand{\dt}{\delta}

\newcommand{\zt}{\zeta}
\newcommand{\et}{\eta}

\newcommand{\te}{\theta}
\newcommand{\ld}{\lambda}

\newcommand{\ph}{\varphi}
\newcommand{\ps}{\psi}
\newcommand{\rh}{\rho}
\newcommand{\om}{\omega}
\newcommand{\ta}{\tau}

\newcommand{\Gm}{\Gamma}

\newcommand{\Q}{{\mathbf{Q}}}
\newcommand{\Z}{{\mathbf{Z}}}
\newcommand{\R}{{\mathbf{R}}}
\newcommand{\C}{{\mathbf{C}}}
\newcommand{\N}{{\mathbf{N}}}

\newcommand{\Diff}{{\mathrm{Diff}}}

\newcommand{\id}{{\mathrm{id}}}

\newcommand{\Ci}{C^{\infty}}

\newcommand{\Aff}{{\mathrm{Aff}}}
\newcommand{\Ker}{{\mathrm{Ker}}}

\newcommand{\andeqn}{\,\,\,\,\,\, {\mbox{and}} \,\,\,\,\,\,}
\newcommand{\QED}{\rule{0.4em}{2ex}}

\newcommand{\csmf}{compact smooth manifold}
\newcommand{\ca}{C*-algebra}
\newcommand{\ct}{continuous}
\newcommand{\pj}{projection}
\newcommand{\mf}{manifold}

\newcommand{\diff}{diffeomorphism}

\newcommand{\hm}{homomorphism}

\newcommand{\ifo}{if and only if}
\newcommand{\rsha}{recursive subhomogeneous algebra}

\newcommand{\hme}{homeomorphism}
\newcommand{\mh}{minimal homeomorphism}
\newcommand{\tgca}{transformation group \ca}
\newcommand{\cp}{crossed product}
\newcommand{\cms}{compact metric space}
\newcommand{\chs}{compact Hausdorff space}
\newcommand{\cfn}{continuous function}

\newcommand{\rsz}[1]{\raisebox{0ex}[0.8ex][0.8ex]{$#1$}}

\title[Examples]{Examples of different minimal diffeomorphisms
giving the same C*-algebras}

\author{N.\  Christopher Phillips}

\address{Department of Mathematics, University  of Oregon,
       Eugene OR 97403-1222, USA.}

\email[]{ncp@darkwing.uoregon.edu}

\thanks{Research
      partially supported by NSF grants DMS 9706850, DMS 0070776,
      DMS 0302401,
      and by the MSRI}

\date{12 August 2004}

\begin{document}

\setcounter{section}{-1}

\begin{abstract}
We give examples of \md s of compact connected \mf s
which are not topologically orbit equivalent,
but whose \tgca s are isomorphic.
The examples show that the following properties of a \md\  are
not invariants of the \tgca:
having topologically quasidiscrete spectrum;
the action on singular cohomology (when the \mf s are diffeomorphic);
the homotopy type of the \mf\  %
(when the \mf s have the same dimension);
and the dimension of the \mf.

These examples also give examples of nonconjugate isomorphic
Cartan subalgebras, and of nonisomorphic Cartan subalgebras,
of simple separable nuclear unital \ca s with tracial rank zero
and satisfying the Universal Coefficient Theorem.
\end{abstract}

\maketitle

\section{Introduction}\label{Sec:Intro}

\indent
The purpose of this paper is to give examples of
distinct \md s of compact connected \mf s
whose \tgca s are isomorphic.
This has become possible because of recent results
which show that, in the real rank zero case,
the \tgca s of \md s are classifiable
in the sense of the Elliott program~\cite{El5}.

For \mh s of the Cantor set,
a remarkable theorem of Giordano, Putnam, and Skau
(Theorem~2.1 of~\cite{GPS}) asserts that
isomorphism of the \tgca s is equivalent to strong orbit equivalence
of the \hme s.
On the circle, it is a consequence of the classification of
the irrational rotation algebras that two \mh s have
isomorphic crossed products \ifo\  they are flip conjugate.
(See the beginning of Section~1 of~\cite{Ph} for details.)

Until now, nothing definite has been known about the
problem on \cms s other than the Cantor set and the circle.
(There have been suggestive results.
See Section~1 of~\cite{Ph} for an extensive discussion.)
In this paper, we show that
strong orbit equivalence is much too strong a relation
to correspond to isomorphism of the \ca s.
We give four kinds of examples of pairs of \md s
of compact connected smooth \mf s which are not even
topologically orbit equivalent,
but which have isomorphic \ca s.
These examples show that, for \md s of the same \mf,
the action of the \hme\  on singular cohomology and
the property of having topologically quasidiscrete spectrum
are not invariants of the \ca;
for \md s of \mf s of the same dimension,
the homotopy type of the \mf\  is not an invariant of the \ca;
and for arbitrary \md s,
the dimension of the \mf\  is not an invariant of the \ca.

For \mh s of a connected \cms, strong orbit equivalence
(even topological orbit equivalence)
turns out to imply flip conjugacy.
(See Theorem~3.1 and Remark~3.4 of~\cite{BT},
or Proposition~5.5 of~\cite{LP}.)
Thus, our examples are equivalently examples of \md s
which have isomorphic \tgca s but which are not flip conjugate.

Our examples have another important consequence.
They provide examples
of simple separable nuclear unital \ca s,
which even have tracial rank zero
in the sense of \cite{LnTAF},~\cite{LnTTR},
with nonconjugate, even nonisomorphic, Cartan subalgebras.
Cartan subalgebras (in the von Neumann algebra sense)
of the hyperfinite type~II$_1$ factor are unique up to
automorphisms of the factor (\cite{CFW}),
although this is not true in a general type~II$_1$ factor~(\cite{CJ}).
Two slightly different definitions of Cartan subalgebras
in \ca s have been introduced, in~\cite{Rn} and~\cite{Kj2}.
Corollary~1.16 in Chapter~III of~\cite{Rn} gives conditions
under which two Cartan subalgebras of a \ca\  must be
conjugate by an automorphism.
(The conditions are restrictive: among other things,
both subalgebras must be~AF.
It isn't part of the hypotheses,
but one sees from the proof that $A$ must also be~AF.)

The existence of nonconjugate Cartan subalgebras is not new:
see the end of the introduction to~\cite{Kj1}
for nonisomorphic Cartan subalgebras,
and see Theorem~\ref{UctNonconjCartan} below,
which is immediate from results already in the literature,
for isomorphic but nonconjugate Cartan subalgebras.
Our examples give both simple nuclear \ca s with
new kinds of isomorphic but nonconjugate Cartan subalgebras,
and simple nuclear \ca s with many Cartan subalgebras
with quite different structure,
being algebras of \cfn s on compact spaces of different dimensions.

All our examples are smooth.
Since the announcement of our results (Section~2 of~\cite{Ph}),
an easier proof
of the required classification results has been found~\cite{LhP},
which also does not require smoothness.
This has made possible
the last three examples in Section~5 of~\cite{LhP},
which give further examples
of pairs of \mh s with isomorphic \tgca s
but which are not topologically orbit equivalent.
Two of these examples, both discussed in Section~1 of~\cite{Ph}
but with the isomorphisms only conjectured,
involve zero and one dimensional \cms s which are highly disconnected.
The third, based on~\cite{IO},
involves connected one dimensional spaces
which are not locally connected.

In this paper, we could omit some of the work by ignoring smoothness.
We prefer to keep smoothness,
because we believe that the problem of finding conditions on \md s
for the isomorphism of smooth crossed products is important.
In particular, this problem might have an answer
quite different from the problem for C*~crossed products.
See Section~3 of~\cite{Ph} for more.
We only note here that for the smooth crossed product algebra to
exist, a kind of temperedness condition must be satisfied.

To summarize,
besides providing examples of different \mh s with isomorphic \tgca s,
there are three ways one can think of the results:
\begin{itemize}
\item
Properties of \mh s which are not invariants of the \tgca.
\item
Examples of nonconjugate,
sometimes nonisomorphic Cartan subalgebras.
\item
Examples on which to test the question of whether smooth
crossed products preserve more information than C*~crossed products.
\end{itemize}
We may also think of the paper as an application of H.~Lin's
classification theorem~\cite{LnCls}.

This paper is organized as follows.
In Section~1, we state in a convenient form some general results
that will be repeatedly used.
Each of the remaining four sections is devoted to one example.

I am grateful to Qing Lin for suggesting a nice generalization
of the original version of the second example.

\section{Preliminaries}\label{Sec:Pre}

In this section, we present for convenience some results
which will be used repeatedly in the rest of the paper.
They are the computation of the ordered K-theory of crossed
products, using the Pimsner-Voiculescu exact sequence
and Exel's rotation numbers,
classification for crossed products by \md s,
some basic facts about flow equivalence,
and some basic facts about Cartan subalgebras of crossed products.
Nothing here is really new.

For reference, we state here the Pimsner-Voiculescu exact sequence
\cite{PV} for the special case
of a crossed product of a compact space by a \hme.
This is what we use to compute unordered K-theory of \cp s.

\begin{thm}\label{PVSeq}
Let $X$ be a \chs, and let $h \colon X \to X$ be a \hme.
Then there is a natural $6$~term exact sequence
\[
\begin{CD}
K^0 ( X ) @>{\id - h^*}>> K^0 ( X ) @>>> K_0 ( C^* (\Z, X, h) )\\
@A{\mathrm{exp}}AA & &  @VV{\partial}V\\
K_1 ( C^* (\Z, X, h) )  @<<< K^1 ( X ) @<<{\id - h^*}< K^1 ( X )
\end{CD}
\]
The maps $K^j (X) \to K_j ( C^* (\Z, X, h) )$ are from the
inclusion $C (X) \to C^* (\Z, X, h).$
\end{thm}

To get the sequence in this form,
take the action on $C (X)$ to be $\af (f) = f \circ h^{-1},$
and identify $K_* (C (X))$ with $K^* (X)$ and
$\af_* \colon K_* ( C (X) ) \to K_* ( C (X) )$ with
$\left( h^{-1} \right)^* \colon K^* (X) \to K^* (X).$

Since we will make frequent use of it, we recall the machinery of
\cite{Ex} for the computation of the map on $K_0$ of a crossed product
by a \hme\  determined by an invariant measure.

\begin{dfn}\label{RDfn}
Let $X$ be a connected \cms, and let $h \colon X \to X$ be a \hme.
Let $A = C^* (\Z, X, h)$ be the \tgca.
Let $\mu$ be an $h$-invariant Borel probability measure on $X.$
Let
\[
K^1 (X)^h = \{ \et \in K^1 (X) \colon h^* (\et) = \et \}.
\]
Define $\rh_h^{\mu} \colon K^1 (X)^h \to S^1$ as follows.
Given $u \in U (C (X, M_n))$ such that $[u \circ h^{-1}] = [u],$
set $z (x) = \det (u (x)),$ find a \cfn\  $a \colon X \to \R$ such that
$z \left( h^{-1} (x) \right)^* z (x) = \exp (2 \pi i a (x)),$
and define
\[
\rh_h^{\mu} ( [u] ) = \exp \left( 2 \pi i \int_X a \, d \mu \right).
\]
The number $\rh_h^{\mu} ( [u] )$ is called the {\emph{rotation number}}
of $[u]$ with respect to $h_n$ and $\mu.$
\end{dfn}

The main result that we use is the following theorem.
It is obtained by combining
Definition~IV.1 and Theorems~V.12 and~VI.11 of~\cite{Ex},
using the definitions and properties of the maps
$\Det_*$ and $R_h^{\mu}$ of Theorem~VI.11 of~\cite{Ex},
which are contained in Section~VI of~\cite{Ex}.
(By comparing Definition~VI.2 with the proof of Proposition~VI.10
in~\cite{Ex}, one sees that the automorphism of $C ( X )$
given by $h$ really is $\af (f) = f \circ h^{-1}.$)

\begin{thm}\label{RThm}
(\cite{Ex})
Let $X,$ $h,$ $A,$ and $\mu$ be as in Definition~\ref{RDfn}.
The function $\rh_h^{\mu} \colon K^1 (X)^h \to S^1$
of Definition~\ref{RDfn}
is a well defined group \hm.
Moreover, if $\partial \colon K_0 ( A) \to K^1 (X)$
is the connecting \hm\  in the Pimsner-Voiculescu exact sequence,
$\ta \colon A \to \C$ is the tracial state on $A$ determined by $\mu,$
and $\et \in K_0 (A),$ then $\partial (\et) \in K^1 (X)^h$ and
$\exp ( 2 \pi i \ta_* (\et)) = \rh_h^{\mu} \circ \partial (\et).$
\end{thm}

Note that by combining the Pimsner-Voiculescu exact sequence,
Exel's results, and Theorem~4.5(1) of~\cite{Ph7}, we have machinery
for the complete computation of the Elliott invariant of the
\tgca\  of a \mh\  of a finite dimensional connected \cms.

The following result is a special case of results of~\cite{LP2},
and is what we use to establish isomorphisms of crossed products.

\begin{thm}\label{Isom}
For $k = 1, \, 2$ let $M_k$ be a \csmf, and let
$h_k \colon M_k \to M_k$ be a uniquely ergodic \md.
Let $\ta_k$ be the unique tracial state on $C^* (\Z, M_k, h_k).$
Assume that there are isomorphisms
\[
f_j \colon K_j ( C^* (\Z, M_1, h_1) ) \to K_j ( C^* (\Z, M_2, h_2) )
\]
such that $f_0 ([1]) = [1]$ and
$(\ta_2)_* \circ f_0 = (\ta_1)_*$ as maps from
$K_0 ( C^* (\Z, M_1, h_1) )$ to $\R.$
Suppose further that $(\ta_1)_*$ has dense range.
Then $C^* (\Z, M_1, h_1) \cong C^* (\Z, M_2, h_2).$
\end{thm}

\begin{proof}
By the main result of~\cite{LP2}
(also see the survey article~\cite{LP3}),
the \cp s are direct limits, with no dimension
growth, of \rsha s.
Theorem~2.3 of~\cite{Ph7} therefore shows that the order on
$K_0$ is determined by the tracial states:
$\et \in K_0 ( C^* (\Z, M_j, h_j) )$ is positive \ifo\  either
$\et = 0$ or $(\ta_j)_* (\et) > 0.$
So $f_0$ is an order isomorphism.
The rest of the hypotheses now imply that the Elliott invariants
of the two \cp s are isomorphic.
Since there is a unique tracial state, the hypothesis that
$(\ta_1)_*$ have dense range implies that,
for $A = C^* (\Z, M_1, h_1),$
the canonical map $\rh_A \colon K_0 (A) \to \Aff (T (A))$
has image $\rh_A ( K_0 (A) )$ dense in $\Aff (T (A)).$
Hence the same is true for $A = C^* (\Z, M_2, h_2).$
Therefore the crossed products have tracial rank zero
by Theorem~4.4 of~\cite{Ph8}.
As in the proof of Theorem~4.5 of~\cite{Ph8},
the classification theorem, Theorem~5.2 of~\cite{LnCls}, applies,
and $C^* (\Z, M_1, h_1) \cong C^* (\Z, M_2, h_2).$
\end{proof}

We can also use Corollary~4.8 of~\cite{LhP}
in place of~\cite{LP2} and~\cite{Ph8}.

We now turn to flow equivalence, a weaker relation than conjugacy.
Since the relevant relation for our purposes is actually
flip conjugacy, we will actually introduce and use
flip flow equivalence.
Flip flow equivalence
does not imply isomorphism of the crossed products,
only Morita equivalence
(Lemma~1.2 of~\cite{Pc1}; Situation~8 of~\cite{Rf}).
Indeed, flip flow equivalence
seems to play the same role for Morita equivalence
that flip conjugacy does for isomorphism.
(Theorem~2.6 of~\cite{GPS}, based on Definition~1.8~\cite{GPS},
is at least suggestive, although in a different context.)
Thus, we give flip flow equivalence less emphasis.
However, it seems inappropriate to ignore it entirely.

The following definitions (except flip flow equivalence) can be found,
for example, at the beginning of Section~1 of~\cite{Pc2} and
in Definition~1.1 and the following discussion in~\cite{Pc1}.

\begin{dfn}\label{CSec}
A {\emph{cross section}} of an action of $\R$ on a \cms\  $X$ is a
closed subset $K \subset X$ such that the map from $\R \times K$ to $X,$
given by $(t, x) \mapsto t x,$ is a surjective local \hme.

We always take the cross section $K$ to be equipped with the
\hme\  $h \colon K \to K$ defined as follows.
For $x \in K,$ let $\ld_K (x) \in \R$ be given by
\[
\ld_K (x) = \inf ( \{ t \in (0, \infty) \colon t \cdot x \in K \} ),
\]
and take $h (x) = \ld_K (x) \cdot x.$
\end{dfn}

Theorem~1 of~\cite{Sz} gives some equivalent conditions for $K$ to
be a cross section.

\begin{dfn}\label{FlEq}
Homeomorphisms $h_1 \colon X_1 \to X_1$ and $h_2 \colon X_2 \to X_2$ are
{\emph{flow equivalent}} if both can be obtained as cross sections
to the same flow.
They are {\emph{flip flow equivalent}} if $h_1$ is flow equivalent to
$h_2$ or to $h_2^{-1}.$
\end{dfn}

One can simplify a little (still following Section~1 of~\cite{Pc2}):

\begin{dfn}\label{Susp}
Let $X$ be a topological space, and let $h \colon X \to X$ be a \hme.
The {\emph{suspension}} of $(X, h)$ is the flow
(action of $\R$ on a topological space) constructed as follows.
On the space $X \times \R,$ we define an action of $\R$
by $s \cdot (x, t) = (x, \, s + t),$ and an action of $\Z$ by taking
the generator to be the \hme\  $(x, t) \mapsto (h (x), \, t - 1).$
These actions commute, and therefore the action of $\R$ on the
space $Y = (X \times \R) / \Z$ is well defined and \ct.
We define the suspension of $(X, h)$
to be $Y$ equipped with this action of $\R.$
\end{dfn}

In~\cite{Pc2}, this flow is called the mapping torus flow.

\begin{lem}\label{FlEqCond}
Homeomorphisms $h_1 \colon X_1 \to X_1$ and $h_2 \colon X_2 \to X_2$ are
flow equivalent \ifo\  $(X_2, h_2)$ can be obtained as a cross
section of the suspension of $(X_1, h_1).$
\end{lem}

Finally, we state some results on Cartan subalgebras.
The following version of a Cartan subalgebra is taken from
Definitions~1 and~3 in Section~1 of~\cite{Kj2}.

\begin{dfn}\label{CSDiagDfn}
Let $A$ be a unital \ca.
A {\emph{normalizer}} of a C*-subalgebra $B$ of $A$ is an element
$a \in A$ such that $a B a^* \subset B$ and $a^* B a \subset B.$
It is {\emph{free}} if $a^2 = 0.$
A {\emph{diagonal}} in $A$ is a commutative C*-subalgebra
$B$ of $A$ which contains $1_A$ and such that
there is a faithful conditional expectation $P \colon A \to B$
whose kernel is the closed linear span of the free normalizers of $B.$
\end{dfn}

By Proposition~4 in Section~1 of~\cite{Kj2},
the definition implies that
a diagonal is a maximal abelian subalgebra.
In the next theorem, we are interested in \mh s of infinite \cms s,
but the proof is the same in greater generality.
One should be able to obtain the result as a corollary
of results in~\cite{Kj2},
but a direct proof seems more illuminating.

\begin{thm}\label{CXIsDiag}
Let $X$ be a \chs\  with a free action of a discrete group $\Gm.$
Then $C (X)$ is a diagonal in $C^* (\Gm, X)$ in the sense
of Definition~\ref{CSDiagDfn}.
\end{thm}

\begin{proof}
The conditional expectation is the standard one:
letting $u_{\gm} \in C^* (\Gm, X)$ be the standard unitary
corresponding to $\gm \in \Gm,$
if all but finitely $f_{\gm} \in C (X)$ are zero, then
$P \left( \sum_{\gm \in \Gm} f_{\gm} u_{\gm} \right) = f_1.$
The only part requiring proof is that its
kernel is the closed linear span of the free normalizers of $C (X).$

First, let $a \in C^* (\Gm, X)$ be a free normalizer;
we prove $P (a) = 0.$
The definition of a normalizer implies in particular that
$a^* a \in C (X).$
Using the properties of conditional expectations at the
second step and $a^2 = 0$ at the third, we get
\[
[a P (a)]^* [a P (a)]
 = P (a)^* a^* a P (a)
 = P (a^*) P ((a^* a) a)
 = 0,
\]
whence $a P (a) = 0.$
Using $P (a) \in C (X),$ we now get
$P (a)^2 = P (a P (a)) = 0.$
Since $P (a)$ is an element of a commutative \ca,
this implies that $P (a) = 0.$

We finish the proof by showing that the closed linear span $L$
of the free normalizers contains ${\mathrm{Ker}} (P).$
Since finite sums of the form $\sum_{\gm \in \Gm} f_{\gm} u_{\gm}$
are dense in $C^* (\Gm, X),$
and since
\[
{\mathrm{Ker}} (P) = \{ a - P (a) \colon a \in C^* (\Gm, X) \},
\]
it suffices to show that $f u_{\gm} \in L$ for every
$f \in C (X)$ and $\gm \in \Gm \setminus \{ 1 \}.$
Since $X$ is compact and the group element $\gm$ has no fixed points,
there is a finite cover of $X$ consisting of open sets
$U \subset X$ such that $\gm U \cap U = \varnothing.$
Choose a partition of unity $(g_1, g_2, \ldots, g_n)$
subordinate to such a cover.
Then $f u_{\gm} = \sum_{k = 1}^n f g_k u_{\gm},$
we have
\[
(f g_k u_{\gm})^2
  = f (f \circ \gm^{-1}) g_k (g_k \circ \gm^{-1}) u_{\gm^2}
  = 0
\]
because
${\mathrm{supp}} (g_k) \cap {\mathrm{supp}} (g_k \circ \gm^{-1})
   = \varnothing,$
and it is trivial that $f g_k u_{\gm}$ normalizes $C (X).$
Thus $f u_{\gm} \in L.$
\end{proof}

The same result also holds
for Renault's definition of a Cartan subalgebra.

\begin{thm}\label{CXIsCartan}
Let $X$ be an infinite \cms, and let $h \colon X \to X$ be a \mh.
Then $C (X)$ is a Cartan subalgebra of $C^* (\Z, X, h)$ in the sense
of Definition~4.13 in Chapter~II of~\cite{Rn}.
\end{thm}

\begin{proof}
This is immediate from Proposition~4.14 in Chapter~II of~\cite{Rn}.
\end{proof}

For the application to our examples,
we combine the above with a result of Tomiyama.

\begin{prp}\label{NonconjCartan}
Let $h_1 \colon X_1 \to X_1$ and $h_2 \colon X_2 \to X_2$
be \mh s of infinite \cms s
such that $h_1$ and $h_2$ are not flip conjugate.
Suppose there is an isomorphism
$\ph \colon C^* (\Z, X_1, h_1) \to C^* (\Z, X_2, h_2).$
Then $\ph (C (X_1))$ and $C (X_2)$ are diagonals in $C^* (\Z, X_2, h_2),$
and are Cartan subalgebras in the sense of Renault,
which are not conjugate by any automorphism of $C^* (\Z, X_2, h_2).$
\end{prp}

\begin{proof}
That $\ph (C (X_1))$ and $C (X_2)$
are diagonals is Theorem~\ref{CXIsDiag},
and that they are Cartan subalgebras is Theorem~\ref{CXIsCartan}.
If there is an automorphism $\ps$ of $C^* (\Z, X_2, h_2)$
such that $\ps (\ph (C (X_1))) = C (X_2),$
then $\ps \circ \ph \colon C^* (\Z, X_1, h_1) \to C^* (\Z, X_2, h_2)$
is an isomorphism sending $C (X_1)$ to $C (X_2),$
so the corollary at the end of~\cite{Tm}
implies that $h_1$ and $h_2$ are flip conjugate.
\end{proof}

We can use results already in the literature to give
examples of a simple unital AT~algebras with real rank zero
which have uncountably many isomorphic but nonconjugate
Cartan subalgebras.

\begin{thm}\label{UctNonconjCartan}
Let $A$ be a simple unital AT~algebra with real rank zero
and with $K_1 (A) \cong \Z.$
Then there exist uncountably many diagonals (Cartan subalgebras)
$B_t \subset A,$ for $t \in [0, \infty],$
such that $B_s \cong B_t$ for all $s$ and $t,$
but such that for $s \neq t$ there is no automorphism $\ph$ of $A$
with $\ph (B_s) = B_t.$
\end{thm}

\begin{proof}
It is not possible to have $K_0 (A) \cong \Z.$
So Theorem~1.15 of~\cite{GPS} provides
a \mh\  $h$ of the Cantor set $X$ such that $A \cong C^* (\Z, X, h).$
Use Theorem~2.3 of~\cite{BH}, Theorem~7.1 of~\cite{Sg2},
and Theorem~6.1 of~\cite{Sg1} to choose,
for every $t \in [0, \infty],$
a \mh\  $h_t$ of $X$ which is strong orbit equivalent to $h$
and which has topological entropy equal to $t.$
Theorem~2.1 of~\cite{GPS} provides isomorphisms
$\ps_t \colon C^* (\Z, X, h_t) \to A.$
Set $B_t = \ps_t (C (X)).$
For $s \neq t,$ the \hme s $h_s$ and $h_t$ are not
flip conjugate because they have different topological entropies.
So Proposition~\ref{NonconjCartan} implies that
$B_s$ and $B_t$ are diagonals (Cartan subalgebras)
which are not conjugate by an automorphism of $A.$
\end{proof}

Applying the isomorphism result of Example~5.8 of~\cite{LhP}
in the same way,
we see that each such algebra also has at least one
diagonal (Cartan subalgebra) which is not isomorphic to
the ones in the theorem.

\section{First example: Furstenberg
        transformations on $(S^1)^2$}\label{Sec:FurstT2}

We give two Furstenberg transformations on the $2$-torus
$S^1 \times S^1$ with isomorphic \tgca s, such that one has
topologically quasidiscrete spectrum and the other does not.
In terms of the three additional ways to think of the results:
\begin{itemize}
\item
Having topologically quasidiscrete spectrum or not is
not an invariant of the \tgca.
\item
The \ca\  in the example
has two isomorphic diagonals (Cartan subalgebras)
which are the algebras of \cfn s on a connected space and
which are not conjugate by an automorphism of the algebra.
\item
Both diffeomorphisms are tempered in the sense of
Definition~3.1 of~\cite{Ph} (see Example~3.7 of~\cite{Ph}),
and the actions are very similar,
so isomorphism of the smooth crossed products
is a very interesting problem.
\end{itemize}

\begin{exa}\label{Rouhani}
Let $\te \in [0, 1] \setminus \Q$ be an irrational number,
and let $r \colon S^1 \to \R$ be a smooth function,
both chosen according to the proof of Lemma~2.3 of~\cite{Rh}.
(See below for details.)
Define $h_1, \, h_2 \colon S^1 \times S^1 \to S^1 \times S^1$ by
\[
h_1 (\zt_1, \zt_2)
 = \ts{ \left( e^{2 \pi i \te} \zt_1, \, \zt_1 \zt_2 \right)  }
\andeqn
h_2 (\zt_1, \zt_2)
 = \ts{ \left( e^{2 \pi i \te} \zt_1,
   \, e^{2 \pi i r (\zt_1) } \zt_1 \zt_2 \right) }
\]
for $(\zt_1, \zt_2) \in S^1 \times S^1.$
(The only difference is the extra factor $\exp (2 \pi i r (\zt_1) )$
in the definition of $h_2.$)

With these choices, it is proved in~\cite{Rh} that $h_1$ has
topologically quasidiscrete spectrum (see Section~1 of~\cite{Rh}
for the definition), while $h_2$ does not.
Therefore $h_1$ is not flip conjugate to $h_2.$
The Elliott invariants of the \tgca s for all such \hm s
are computed in Example~4.9 of~\cite{Ph7}, and in particular
it is shown there that for uniquely ergodic \hme s of this form,
the Elliott invariants are all isomorphic,
and that the unique tracial state $\ta$ has
$\ta_* (K_0 ( C^* (\Z, \, S^1 \times S^1, \, h )))$ dense in $\R.$
Now $h_2$ is uniquely ergodic
(see the proof of Theorem~2.1 of~\cite{Rh}),
and $h_1$ is uniquely ergodic and both $h_1$ and $h_2$ are minimal
(see~\cite{Rh} or Example~4.9 of~\cite{Ph7}; the original source
is Section~2 of~\cite{Fr}).
Therefore, to use Theorem~\ref{Isom} to prove that
$C^* (\Z, \, S^1 \times S^1, \, h_1 ) \cong
         C^* (\Z, \, S^1 \times S^1, \, h_2 ),$
it suffices to verify that $r$ is smooth.

Recall from the proof of Lemma~2.3 of~\cite{Rh} that
$\nu_1 = 1,$ that
$\nu_{k + 1} = 2^{\nu_k} + \nu_k + 1$
for $k \geq 1,$ and that
$n_k = {\mathrm{sgn}} (k) \cdot 2^{\nu_{| k |}}$
for $k \in \Z \setminus \{ 0 \}.$
Further recall that
\[
\te = \sum_{k = 1}^{\infty} 2^{- \nu_k}
\]
and
\[
r (t) = \sum_{k \in \Z \setminus \{ 0 \} } \bt_k e^{2 \pi i n_k t},
\]
with
\[
\bt_k = \frac{1}{ |k| } \left( e^{2 \pi i n_k \te} - 1 \right).
\]
It follows from the proof in~\cite{Rh} that
$| \bt_k | \leq 2 \pi \cdot 2^{- | n_k|} \cdot |k|^{-1}$
for all $k \neq 0.$
To prove that $r$ is smooth, it is enough to prove uniform convergence,
for every $m \geq 0,$ of the series
\[
\sum_{k \in \Z \setminus \{ 0 \} }
      (2 \pi i n_k)^m \bt_k e^{2 \pi i n_k t},
\]
obtained by differentiating the series for $r (t)$
term by term $m$ times.
For $k \geq 1$ the $n_k$ are distinct positive integers, so
\[
\sum_{k \in \Z \setminus \{ 0 \} } \ts{ \left| (2 \pi i n_k)^m \bt_k \right| }
 \leq 2 (2 \pi)^{m + 1} \sum_{k = 1}^{\infty} \frac{n_k^m}{2^{n_k} k}
 \leq 2 (2 \pi)^{m + 1} \sum_{n = 1}^{\infty} \frac{n^m}{2^{n}}
     < \infty.
\]
Since the functions $t \mapsto e^{2 \pi i n_k t}$
have absolute value $1,$
this proves the required uniform convergence.

It now follows from Theorem~\ref{Isom} that
$C^* (\Z, \, S^1 \times S^1, \, h_1 ) \cong
         C^* (\Z, \, S^1 \times S^1, \, h_2 ).$

{}From Proposition~\ref{NonconjCartan},
we see that the $C^* (\Z, \, S^1 \times S^1, \, h_2 )$
has two nonconjugate diagonals (Cartan subalgebras),
both isomorphic to $C (S^1 \times S^1).$
\end{exa}

The example answers the isomorphism question
raised in~\cite{Rh} before Proposition~2.5,
but in the opposite way to what is suggested there.

It is likely that the \diff s $h_1$ and $h_2$ can be shown not to be
flip flow equivalent (Definition~\ref{FlEq}),
by using Theorem~4.1 of~\cite{Pc2}.

One intriguing question arises.

\begin{qst}
Let $h_1$ and $h_2$ be as in Example~\ref{Rouhani}.
Does there exist a minimal \hme\  $h$ of $S^1 \times S^1$ such that
both $h_1$ and $h_2$ are factors of $h$?
\end{qst}

\section{Second example: Affine Furstenberg
        transformations on $(S^1)^3$}\label{Sec:FurstT3}

We give two affine Furstenberg transformations on
$(S^1)^3$
whose \ca s are isomorphic but which are not flip conjugate.
Variations produce an arbitrarily large number of such
transformations.
The isomorphism of the \ca s answers a question raised in
Section~6.1 of the unpublished thesis of R.\  Ji~\cite{Ji}.
We point out here that there is a mistake in~\cite{Ji},
which is identified and corrected in followup work~\cite{RM}.
In particular, we refer to~\cite{RM} for the best current
knowledge of the (unordered) K-theory of \tgca s of
Furstenberg transformations in high dimensions.
The mistake in~\cite{Ji} does not affect the examples in this section.

In terms of the three additional ways to think of the results:
\begin{itemize}
\item
The action of the \hme\  on integral cohomology
not an invariant of the \tgca.
\item
There are \ca s with arbitrarily large finite numbers
of isomorphic diagonals (Cartan subalgebras)
which are the algebras of \cfn s on a connected space and
which are not conjugate by automorphisms of the algebra.
\item
All the diffeomorphisms are tempered in the sense of
Definition~3.1 of~\cite{Ph} (see Example~3.6 of~\cite{Ph}),
and the actions are very similar,
so isomorphism of the smooth crossed products
is a very interesting problem.
\end{itemize}

We begin with a general calculation for affine Furstenberg
transformations on $(S^1)^3.$
It is mostly a special case of computations done in~\cite{Ji};
the only new part is the determination of the order on the K-theory
of the crossed product, rather than merely what the (unique)
tracial state does on K-theory.
A similar calculation for $(S^1)^2$ was done in
Example~4.9 of~\cite{Ph7}.
We give the calculation here because~\cite{Ji} has never been
published, and because there are new features in this case which do
not appear in the calculation of~\cite{Ph7}.
Also, we will need this result for our third example.

\begin{lem}\label{FT3D}
Let $M = (S^1)^3,$ let $m, \, n \in \Z \setminus \{ 0 \},$
and define $h = h_{m, n, \te} \colon M \to M$ by
\[
h (\zt_1, \zt_2, \zt_3)
 = \left( \exp (2 \pi i \te) \zt_1, \, \zt_1^m \zt_2,
              \, \zt_2^n \zt_3 \right)
\]
for $(\zt_1, \zt_2, \zt_3) \in (S^1)^3.$
Then $h$ is minimal and uniquely ergodic.
Moreover, the groups $K_0 (C^* (\Z, M, h))$ and $K_1 (C^* (\Z, M, h))$
are both isomorphic to $\Z^4 \oplus \Z / m \Z \oplus \Z / n \Z.$
The isomorphism of $K_0 (C^* (\Z, M, h))$ with this group can
be chosen in such a way that the unique tracial state $\ta$ induces
the map
\[
\ta_* ( r_1, r_2, r_3, r_4, s_1, s_2) = r_1 + \te r_3
\]
and $K_0 (C^* (\Z, M, h))_+$ is identified with
\[
\{ ( r_1, r_2, r_3, r_4, s_1, s_2)
    \in \Z^4 \oplus \Z / m \Z \oplus \Z / n \Z \colon
     r_1 + r_3 \te > 0 \}
   \cup \{ 0 \}.
\]
\end{lem}

\begin{proof}
That $h$ is minimal and uniquely ergodic follows from Theorem~2.1
in Section~2.3 of~\cite{Fr}.
By this same theorem, the unique ergodic measure is the normalized
Lebesgue measure $\ld$ on $(S^1)^3.$

The K\"{u}nneth Theorem (Theorem~4.1 of~\cite{Sc1}; also see
Corollary 2.7.15 of~\cite{At} for the commutative case,
which suffices here)
shows that
\[
K^* {\ts{ \left( \rsz{ \ts{ (S^1)^3 } } \right) }}
 \cong K^* ( S^1) \otimes K^* ( S^1) \otimes K^* ( S^1).
\]
We identify the two sides of this isomorphism.
Let $1 \in C (S^1)$ be the identity, and let $z \in C (S^1)$
be the canonical unitary $z (\zt) = \zt.$
Then $K^0 {\ts{ \left( \rsz{ \ts{ (S^1)^3 } } \right) }}$
is the free abelian group on generators
\[
\et_1 = [1] \otimes [1] \otimes [1], \,\,\,\,\,\,
\et_2 = [z] \otimes [z] \otimes [1],
\]
\[
\et_3 = [z] \otimes [1] \otimes [z], \andeqn
\et_4 = [1] \otimes [z] \otimes [z],
\]
and $K^1 {\ts{ \left( \rsz{ \ts{ (S^1)^3 } } \right) }}$
is the free abelian group on generators
\[
\gm_1 = [z] \otimes [1] \otimes [1], \,\,\,\,\,\,
\gm_2 = [1] \otimes [z] \otimes [1],
\]
\[
\gm_3 = [1] \otimes [1] \otimes [z], \andeqn
\gm_4 = [z] \otimes [z] \otimes [z].
\]
Moreover, in the graded ring structure on
$K^* {\ts{ \left( \rsz{ \ts{ (S^1)^3 } } \right) }},$
we have
\[
\et_2 = \gm_1 \gm_2, \,\,\,\,\,\,
\et_3 = \gm_1 \gm_3, \,\,\,\,\,\,
\et_4 = \gm_2 \gm_3, \,\,\,\,\,\,
\gm_4 = \gm_1 \gm_2 \gm_3,  \andeqn
\gm_1^2 = \gm_2^2 = \gm_3^2 = 0.
\]
To compute $h^*,$ it therefore suffices to calculate
$h^* (\gm_1),$ $h^* (\gm_2),$ and $h^* (\gm_3).$
We may replace $h$ by the homotopic map
\[
h_0 (\zt_1, \zt_2, \zt_3)
= \left( \zt_1, \, \zt_1^m \zt_2, \, \zt_2^n \zt_3 \right).
\]
Then $h^* (\gm_1)$ is the class of the function
\[
(z \otimes 1 \otimes 1) \circ h_0 = z \otimes 1 \otimes 1,
\]
so $h^* (\gm_1) = \gm_1.$
Similarly,
\[
h^* (\gm_2) = [ (1 \otimes z \otimes 1) \circ h_0 ]
  = [z^m \otimes z \otimes 1] = m \gm_1 + \gm_2
\andeqn
h^* (\gm_3) = n \gm_2 + \gm_3.
\]
It follows that
\[
h^* (\et_2) = \gm_1 (m \gm_1 + \gm_2) = \et_2
\]
(using $\gm_1^2 = 0$), that
\[
h^* (\et_3) = \gm_1 (n \gm_2 + \gm_3 ) = n \et_2 + \et_3,
\]
that
\[
h^* (\et_4) = (m \gm_1 + \gm_2) (n \gm_2 + \gm_3 )
  = m n \et_2 + m \et_3 + \et_4,
\]
and that $h^* (\gm_4) = \gm_4.$
So the matrices of $\id - h^*$ on
$K^0 {\ts{ \left( \rsz{ \ts{ (S^1)^3 } } \right) }}$ and
$K^1 {\ts{ \left( \rsz{ \ts{ (S^1)^3 } } \right) }}$
are given by
\[
0 \oplus
 \left( \begin{array}{ccc}
  0 & -n & -m n \\ 0 & 0 & -m \\ 0 & 0 & 0 \end{array} \right)
\andeqn
\left( \begin{array}{ccc}
  0 & -m & 0 \\ 0 & 0 & -n \\ 0 & 0 & 0 \end{array} \right)
\oplus 0.
\]
Thus
\[
\id - h^* :
 K^0 {\ts{ \left( \rsz{ \ts{ (S^1)^3 } } \right) }}
 \to K^0 {\ts{ \left( \rsz{ \ts{ (S^1)^3 } } \right) }}
\]
has kernel $\Z \et_1 + \Z \et_2$ and cokernel
$\Z \overline{\et}_1 + \Z \overline{\et}_2  + \Z \overline{\et}_3
   + \Z \overline{\et}_4,$
in which the images $\overline{\et}_1$ and $\overline{\et}_4$
of $\et_1$ and $\et_4$ have infinite order, in which
$\overline{\et}_2$ has order $n,$
and in which $\overline{\et}_3$ has order $m.$
Similarly,
\[
\id - h^* :
 K^1 {\ts{ \left( \rsz{ \ts{ (S^1)^3 } } \right) }}
 \to K^1 {\ts{ \left( \rsz{ \ts{ (S^1)^3 } } \right) }}
\]
has kernel $\Z \gm_1 + \Z \gm_4$ and cokernel isomorphic to 
$\Z^2 \oplus \Z / m \Z \oplus \Z / n \Z.$

The exact sequence of Theorem~\ref{PVSeq} therefore breaks apart into
the two exact sequences
\[
0 \longrightarrow [ \Z \overline{\et}_1 + \Z \overline{\et}_4 ]
      \oplus \Z / m \Z \oplus \Z / n \Z
  \longrightarrow K_0 ( C^* (\Z, M, h))
  \stackrel{\partial}{\longrightarrow} \Z \gm_1 + \Z \gm_4
  \longrightarrow 0
\]
and
\[
0 \longrightarrow \Z^2 \oplus \Z / m \Z \oplus \Z / n \Z
  \longrightarrow K_1 ( C^* (\Z, M, h))
  \longrightarrow \Z^2
  \longrightarrow 0.
\]
Both split because $\Z^2$ is free.
Therefore $K_0 ( C^* (\Z, M, h))$ and $K_1 ( C^* (\Z, M, h))$
are both isomorphic to $\Z^4 \oplus \Z / m \Z \oplus \Z / n \Z,$
as claimed.
We are now done with $K_1 ( C^* (\Z, M, h)),$
but it remains to determine
the action of the tracial state and the order
on $K_0 ( C^* (\Z, M, h)).$

Identify $\overline{\et}_1, \, \dots, \, \overline{\et}_4$
with their images in $K_0 ( C^* (\Z, M, h)).$
Let $\ta$ be the unique tracial state on $C^* (\Z, M, h),$ which comes
from normalized Lebesgue measure $\mu$ on $(S^1)^3.$
(See the beginning of the proof.)
Clearly $\ta_* ( \overline{\et}_1 ) = 1.$
Naturality in the K\"{u}nneth formula
for $S^1 \times S^1$ shows that the image under any point
evaluation of
$[z] \otimes [z] \in K_0 ( C (S^1 \times S^1) )$ is zero.
Therefore $[z] \otimes [z]$ can be represented as a difference
$[p] - [q]$ of the classes of two \pj s of the same rank.
(Actually, $[z] \otimes [z]$ is a Bott element, but we do not need
this much.)
The same is therefore true of $\et_2,$ $\et_3,$ and $\et_4,$
whence
$\ta_* ( \overline{\et}_2 ) = \ta_* ( \overline{\et}_3 )
    = \ta_* ( \overline{\et}_4 ) = 0.$

The calculations above show that $K^1 (M)^h = \Z \gm_1 + \Z \gm_4,$
so our next step is to calculate $\rh_h^{\mu} ( \gm_1 )$ and
$\rh_h^{\mu} ( \gm_4 )$ as in Definition~\ref{RDfn}.
We have $\gm_1 = [z \otimes 1 \otimes 1]$ and
\[
(z \otimes 1 \otimes 1) \circ h^{-1}
   = \exp ( - 2 \pi i \te) (z \otimes 1 \otimes 1),
\]
so $\rh_h^{\mu} ( \gm_1 ) = \exp (2 \pi i \te)$ is immediate.
For $\gm_4,$ write $[z] \otimes [z] = [p] - [q]$ as above.
Working in a suitable matrix algebra, we see that $\gm_4$ is represented
by the unitary
\[
[ (1 - p) \otimes 1 + p \otimes z ]
      [(1 - q) \otimes 1 + q \otimes z^{-1} ].
\]
Its determinant is the constant function $1,$
so $\rh_h^{\mu} ( \gm_4 ) = 1.$

Choose $\nu_1^{(0)}, \, \nu_4^{(0)} \in K_0 ( C^* (\Z, M, h))$ such that
$\ts{ \partial \left( \rsz{ \nu_1^{(0)} } \right) } = \gm_1$
and $\ts{ \partial \left( \rsz{ \nu_4^{(0)} } \right) } = \gm_4.$
Theorem~\ref{RThm} implies that
\[
\ts{ \ta_* \left( \rsz{ \nu_1^{(0)} } \right) } \in \te + \Z
\andeqn
\ts{ \ta_* \left( \rsz{ \nu_4^{(0)} } \right) } \in \Z.
\]
Taking $\nu_1 = \nu_1^{(0)} - k \overline{\et}_1$ and
$\nu_4 = \nu_4^{(0)} - l \overline{\et}_1$ for suitable $k$ and $l,$
we get $\ta_* (\nu_1) = \te$ and $\ta_* (\nu_4) = 0.$

We can now identify $K_0 ( C^* (\Z, M, h))$ as
\[
\Z \overline{\et}_1 \oplus \Z \overline{\et}_4 \oplus \Z \nu_1
    \oplus \Z \nu_4 \oplus (\Z / m \Z) \cdot \overline{\et}_2
    \oplus (\Z / n \Z) \cdot \overline{\et}_3,
\]
with
\[
\ta_* ( \overline{\et}_1 ) = 1, \,\,\,\,\,\,
\ta_* ( \nu_1 ) = \te, \andeqn
\ta_* ( \overline{\et}_2 ) = \ta_* ( \overline{\et}_3 )
  = \ta_* ( \overline{\et}_3 ) = \ta_* ( \nu_4 ) = 0.
\]
This is the required formula for $\ta_*,$
and the identification of $K_0 ( C^* (\Z, M, h))_+$ follows from
Theorem~4.5(1) of~\cite{Ph7}.
\end{proof}

In Section~6.1 of~\cite{Ji},
it is stated without proof that if $| m | \neq | n |$
and $\te \in \R \setminus \Q,$
then the \hme s $h_{m, n, \te}$ and $h_{n, m, \te},$
as defined in Lemma~\ref{FT3D},
are not flip conjugate.
(The difference is that $m$ and $n$ are switched.)
We prove this, and in fact a more general statement,
by computing the action on singular cohomology.
The generalization was suggested by Qing Lin.

We start with the following lemma.

\begin{lem}\label{NonSim}
Let
\[
a_1 = \left( \begin{array}{ccc}
  1 & m_1 & r_1 \\ 0 & 1 & n_1 \\ 0 & 0 & 1  \end{array} \right)
\andeqn
a_2 = \left( \begin{array}{ccc}
  1 & m_2 & r_2 \\ 0 & 1 & n_2 \\ 0 & 0 & 1  \end{array} \right)
\]
be integer matrices, with $m_1,$ $n_1,$ $m_2,$ and $n_2$ all nonzero.
Suppose $a_1$ is similar over $\Z$ to $a_2$ or to $a_2^{-1}.$
Then $| m_1 | = | m_2 |$ and $| n_1 | = | n_2 |.$
\end{lem}

\begin{proof}
First assume $a_1$ is similar over $\Z$ to $a_2.$
Set $c_j = a_j - 1.$
Then $c_1$ is similar over $\Z$ to $c_2.$
A calculation shows that
\[
c_j (\Ker (c_j^2) ) = m_j \Z \cdot (1, 0, 0) = m_j \Ker (c_j). 
\]
Therefore $| m_1 | = | m_2 |.$
Another calculation shows that
\[
\Z^3 / c_j^2 ( \Z^3) \cong \Z / m_j n_j \Z \oplus \Z \oplus \Z.
\]
So $| m_1 n_1 | = | m_2 n_2 |.$
Since all four numbers are nonzero and $| m_1 | = | m_2 |,$
it follows that $| n_1 | = | n_2 |.$

Now suppose that $a_1$ is similar over $\Z$ to $a_2^{-1}.$
Since
\[
a_2^{-1} = \left( \begin{array}{ccc}
  1 & - m_2 & m_2 n_2 - r_2 \\ 0 & 1 & - n_2 \\ 0 & 0 & 1  \end{array}
 \right),
\]
the case already considered implies
$| m_1 | = | m_2 |$ and $| n_1 | = | n_2 |$
in this case as well.
\end{proof}

\begin{lem}\label{NonConj}
Let $M = (S^1)^3.$
Let $m_1, n_1, m_2, n_2 \in \Z \setminus \{ 0 \}.$
For $j = 1, 2$
define $h_j = h_{m_j, n_j, \te} \colon M \to M$ by
\[
h_j (\zt_1, \zt_2, \zt_3)
 = \left( \exp (2 \pi i \te) \zt_1, \, \zt_1^{m_j} \zt_2,
              \, \zt_2^{n_j} \zt_3 \right)
\]
for $(\zt_1, \zt_2, \zt_3) \in (S^1)^3.$
If $h_1$ is flip conjugate to $h_2$
then $| m_1 | = | m_2 |$ and $| n_1 | = | n_2 |.$
\end{lem}

\begin{proof}
If $h_{m_1, n_1, \te}$ is conjugate to $h_{m_2, n_2, \te},$
then in particular
there must be an automorphism $b$ of
$H^1 {\ts{ \left( \rsz{ \ts{ (S^1)^3 } }; \, \Z \right) }}$
such that
$b \circ (h_{m_1, n_1, \te})^* \circ b^{-1} = (h_{m_2, n_2, \te})^*.$
If instead $h_{m_1, n_1, \te}$ is conjugate to
$(h_{m_2, n_2, \te})^{-1},$
then we get a similar equation with
$[(h_{m_2, n_2, \te})^*]^{-1}$ on the right.
We show that no such $b$ can exist.

We do the calculation of $(h_{m, n, \te})^*$ in the proof of
Lemma~\ref{FT3D}, but on singular cohomology rather than K-theory,
and using the K\"{u}nneth formula for singular cohomology,
Theorem 5.6.1 of~\cite{Sp}.
We get
$H^1 {\ts{ \left( \rsz{ \ts{ (S^1)^3 } }; \, \Z \right) }}
   \cong \Z^3,$
with a  $\Z$-basis with respect to which
the matrix of $(h_{m, n, \te})^*$ is
\[
(h_{m, n, \te})^* = \left( \begin{array}{ccc}
  1 & m & 0 \\ 0 & 1 & n \\ 0 & 0 & 1  \end{array} \right).
\]
The nonexistence of $b$ now follows from Lemma~\ref{NonSim}.
\end{proof}

\begin{exa}\label{Ex2a}
Fix $\te \in [0, 1] \setminus \Q$ and $m, \, n \in \Z$ with
$0 < m < n.$
Then the two affine Furstenberg transformations on
$(S^1)^3,$ given by
\[
(\zt_1, \zt_2, \zt_3) \mapsto
  \left( \exp (2 \pi i \te) \zt_1, \, \zt_1^m \zt_2,
              \, \zt_2^n \zt_3 \right)
\]
and
\[
(\zt_1, \zt_2, \zt_3) \mapsto
  \left( \exp (2 \pi i \te) \zt_1, \, \zt_1^n \zt_2,
              \, \zt_2^m \zt_3 \right)
\]
(the difference is that $m$ and $n$ have been exchanged),
are not topologically orbit equivalent but have isomorphic
crossed product \ca s.

They are not flip conjugate by Lemma~\ref{NonConj}, so
not topologically orbit equivalent by
Theorem~3.1 and Remark~3.4 of~\cite{BT},
or by Proposition~5.5 of~\cite{LP}.
The Elliott invariants are isomorphic by Lemma~\ref{FT3D},
and both are uniquely ergodic minimal diffeomorphisms
whose tracial states induce maps from $K_0$ to $\R$ with dense
ranges (namely $\Z + \te \Z$).
Therefore the two crossed products are isomorphic by Theorem~\ref{Isom}.

{}From Proposition~\ref{NonconjCartan},
we see that the algebra
has two nonconjugate diagonals (Cartan subalgebras),
both isomorphic to $C ( (S^1)^3).$
\end{exa}

We do not know whether these \hme s are flip flow equivalent
in the sense of Definition~\ref{FlEq}.

\begin{exa}\label{Ex2b}
Let $\te \in [0, 1] \setminus \Q,$
let $r \in \N,$ and let $p_1, p_2, \ldots, p_r$ be distinct primes.
For $0 \leq k \leq r$ set
\[
m_k = p_1 p_2 \cdots p_k \andeqn n_k = p_{k + 1} p_{k + 2} \cdots p_r.
\]
By the same reasoning as in Example~\ref{Ex2a},
the \hme s (in the notation of Lemma~\ref{FT3D})
$h_{m_0, n_0, \te}, \, h_{m_1, n_1, \te}, \ldots, h_{m_r, n_r, \te}$
are pairwise not topologically orbit equivalent,
but, since
\[
\Z / m_k \Z \oplus \Z / n_k \Z \cong \Z / p_1 p_2 \cdots p_r \Z
\]
for all $k,$
all give isomorphic crossed products.
The common crossed product
has $r + 1$ nonconjugate diagonals (Cartan subalgebras),
all isomorphic to $C ( (S^1)^3).$
\end{exa}

\section{Third example: Minimal diffeomorphisms on distinct three
         dimensional manifolds}\label{Sec:Diff3Mfs}

We produce \md s of $S^2 \times S^1$ and of $(S^1)^3$
whose crossed product \ca s are isomorphic.
These are manifolds of the same dimension, but we will see
that it is not possible
for a \md\  on $S^2 \times S^1$ to be flip flow
equivalent to a \md\  on $(S^1)^3,$
let alone flip conjugate.
In terms of the three additional ways to think of the results:
\begin{itemize}
\item
The homotopy type of the space on which the \hme\  acts
is not an invariant of the \tgca,
even if the dimension of the space is fixed.
\item
The \ca\  in the example
has two nonisomorphic diagonals (Cartan subalgebras).
\item
We don't know if the \md\  of $S^2 \times S^1$
can be chosen to be tempered in the sense of
Definition~3.1 of~\cite{Ph},
but it seems reasonable to hope that it can be.
\end{itemize}

\begin{lem}\label{RotNoSC}
Let $X$ be a connected \cms.
Let $u \in M_n (C (X) )$ be unitary.
For $n \in \N$ let $h_n \colon X \to X$ be a \hme\  such that
$h_n^* ([u]) = [u]$ in $K^1 (X),$ and let
$h \colon X \to X$ be a \hme\  such that $h_n \to h$ uniformly.
For each $n \in \N$ let $\mu_n$ be
an $h_n$-invariant Borel probability measure on $X,$
and set $\ld_n = \rh_{h_n}^{\mu_n} ([u]).$
(See Definition~\ref{RDfn}.)
Suppose $\lim_{n \to \infty} \ld_n = \ld.$
Then there exists an $h$-invariant Borel probability measure $\mu$
on $X$ such that $\rh_{h}^{\mu} ([u]) = \ld.$
\end{lem}

\begin{proof}
It is clear from Definition~\ref{RDfn} that the rotation number
is the same for the function $x \mapsto \det (u (x))$ from $X$
to $S^1.$
Thus, we may assume that $u$ itself is a unitary in $C (X).$
Definition~\ref{RDfn} then means that
there are \cfn s $a_n \colon X \to \R$ such that
\[
\ts{ \left( u \circ h_n^{-1} \right)^* } u  = \exp (2 \pi i a_n)
\andeqn
\exp \left( 2 \pi i \int_X a_n \, d \mu_n \right) = \ld_n.
\]

Uniform convergence of $h_n$ to $h$ implies that
$\limi{n} \left\| f \circ h_n - f \circ h \right\| = 0$
for all $f \in C (X).$
Furthermore, we claim that
$\limi{n} \left\| f \circ h_n^{-1} - f \circ h^{-1} \right\| = 0$
for all $f \in C (X).$
Indeed, with $g = f \circ h^{-1},$ we get
\[
\ts{ \left\| f \circ h_n^{-1} - f \circ h^{-1} \right\| }
 = \ts{ \left\| g \circ h \circ h_n^{-1} - g \right\| }
 = \ts{ \left\| \left( g \circ h - g \circ h_n \right)
         \circ h_n^{-1} \right\| }
 \to 0.
\]

The set of Borel probability measures on $X$ is a weak* compact
subset of the dual of $C (X),$ by Alaoglu's Theorem,
and it is weak* metrizable because
it is bounded and $C (X)$ is separable.
By passing to a subsequence, we may therefore assume that
there is a Borel probability measure $\mu$ on $X$ such that
$\mu_n \to \mu$ weak*.

We claim that $\mu$ is $h$-invariant, and we prove this by showing
that
\[
\int_X (f \circ h) \, d \mu = \int_X f \, d \mu
\]
for all $f \in C (X).$
So let $f \in C (X).$
We estimate:
\begin{align*}
& \left| \int_X (f \circ h) \, d \mu - \int_X f \, d \mu \right|  \\
& \hspace*{2em} \mbox{} \leq
 \left| \int_X (f \circ h) \, d \mu
               - \int_X (f \circ h) \, d \mu_n \right|
     + \left| \int_X (f \circ h) \, d \mu_n
                   - \int_X (f \circ h_n) \, d \mu_n \right|
           \\
& \hspace{4em} \mbox{} +
      \left| \int_X
             (f \circ h_n) \, d \mu_n - \int_X f \, d \mu_n \right|
   +  \left| \int_X f \, d \mu_n - \int_X f \, d \mu \right|.
\end{align*}
The third term is zero because $\mu_n$ is $h_n$-invariant,
the first and last terms converge to zero because $\mu_n \to \mu$ weak*,
and the second term converges to zero because
$\limi{n} \left\| f \circ h_n - f \circ h \right\| = 0.$
This proves the claim.

Now set
\[
w_n (x) = u \ts{ \left( \rsz{ h_n^{-1} (x) } \right)^* } u (x)
\andeqn
w (x) = u \ts{ \left( \rsz{ h^{-1} (x) } \right)^* } u (x)
  = \limi{n} w_n (x)
\]
(uniform convergence).
We know that $w_n$ is homotopically trivial in $U (C (X)),$
so $w$ is also.
Therefore there is a \cfn\  $b \colon X \to \R$ such that
$w = \exp (2 \pi i b).$

Fix $x_0 \in X.$
Then $\exp (2 \pi i a_n (x_0)) \to \exp (2 \pi i b (x_0)).$
Therefore there are integers $k_n$ such that
$a_n (x_0) + k_n \to b (x_0).$
Define $b_n = a_n + k_n \in C (X).$
We claim that $\| b_n - b \| \to 0.$

To prove the claim, first notice that if $\zt_1, \, \zt_2 \in S^1$
with $| \zt_1 - \zt_2 | < 1,$ then the length of the arc from
$\zt_1$ to $\zt_2$ is less than
\[
| {\mathrm{Re}} (\zt_1) - {\mathrm{Re}} (\zt_2) |
  + | {\mathrm{Im}} (\zt_1) - {\mathrm{Im}} (\zt_2) |
 \leq 2 | \zt_1 - \zt_2 |.
\]
If therefore $\exp (2 \pi i \af_j) = \zt_j$ for $j = 1, \, 2,$
then there is a unique integer $k$ such that
$| \af_1 + k - \af_2 | < \frac{1}{\pi}$; in fact, one gets
$| \af_1 + k - \af_2 | \leq \frac{1}{\pi} | \zt_1 - \zt_2 |.$
Now, dropping initial terms of the sequence, we may assume that
\[
\| w_n - w \| < 1 \andeqn
  | b_n (x_0) - b (x_0) | < \ts{ \frac{1}{\pi} }
\]
for all $n.$
For any $x \in X$ and $n \in \N,$ let $l_n (x)$ be the unique
integer such that
$| b_n (x) + l_n (x) - b (x) | < \frac{1}{\pi}.$
Since
\[
| b_n (x) + l_n (x) - b (x) |
  \leq \ts{ \frac{1}{\pi} } | w_n (x) - w (x) |
  \leq \ts{ \frac{1}{\pi} } \| w_n - w (x) \|
\]
by the above, we see that $b_n + l_n \to b$ uniformly.
Further, for each $n$ and as $l$ runs through $\Z,$ the sets
\[
\{ x \in X \colon l_n (x) = l \}
  = \ts{ \left\{ x \in X \colon | b_n (x) + l - b (x) |
          < \ts{ \frac{1}{\pi} }  \right\} }
\]
are disjoint, open, and cover $X.$
Because $X$ is connected, only one can be nonempty, necessarily the one
for $l = 0.$
So $\| b_n - b \| \to 0,$ and the claim is proved.

We now claim that
\[
\limi{n} \int_X b_n \, d \mu_n = \int_X b \, d \mu.
\]
Using the weak* convergence $\mu_n \to \mu$ on the second
term at the second step, we have
\[
\left| \int_X b_n \, d \mu_n - \int_X b \, d \mu \right|
  \leq \left| \int_X (b_n - b) \, d \mu_n \right|
    + \left| \int_X b \, d \mu_n - \int_X b \, d \mu \right|
 \to 0
\]
as $n \to \infty.$
This proves the claim.

Now
\[
\ld = \limi{n} \exp \left( 2 \pi i \int_X b_n \, d \mu_n \right)
  = \exp \left( 2 \pi i \int_X b \, d \mu \right),
\]
which by definition is $\rh_h^{\mu} ( [u] ).$
\end{proof}

\begin{lem}\label{hty}
Let $M$ be a connected compact manifold such that $H^1 (M; \Z) = 0.$
Then the group $[M \times S^1, \, S^1]$ of homotopy classes
(with operation given by pointwise multiplication on the codomain)
is isomorphic to $\Z,$ and is generated by the class of the
function $u (x, \zt) = \zt.$
\end{lem}

\begin{proof}
First, note that $M \times S^1,$ being a compact smooth \mf,
is a finite CW~complex.
(See Theorem~3.5 and Corollary~6.7 of~\cite{Ml}.)
Therefore Theorem 8.1.8 and 8.1.1 of~\cite{Sp}
(with $\pi^Y$ being defined in Section~1.3 of~\cite{Sp}) show that
$[M \times S^1, \, S^1] \cong H^1 (M \times S^1; \, \Z).$
We compute $H^1 (M \times S^1; \, \Z)$ using the K\"{u}nneth formula
for singular cohomology, Theorem 5.6.1 of~\cite{Sp}.
This is allowed because $\Z$ is finitely generated and
$H^* (S^1; \Z)$ has finite type.
Since $H^* (S^1; \Z)$ is free and $H^1 (M; \Z) = 0,$ the result is just
\[
H^1 (M \times S^1; \, \Z) \cong H^0 (M; \Z) \otimes H^1 (S^1; \Z)
  \cong H^0 (M; \Z) \cong \Z.
\]
One checks, using the formulas for the maps in~\cite{Sp},
that $u$ corresponds to a generator of $H^1 (M \times S^1; \, \Z).$
\end{proof}

\begin{prp}\label{Exist}
Let $M$ be a connected compact manifold such that $H^1 (M; \Z) = 0.$
Define a unitary $u \in C (M \times S^1)$ by $u (x, \zt) = \zt.$
Then there exists a uniquely ergodic minimal diffeomorphism $h$ of
$M \times S^1$ which is homotopic to the identity map and such
that, with $\mu$ being the unique invariant Borel probability measure,
$\rh_h^{\mu} ([u]) \in S^1$ is equal to $\exp (2 \pi i \te)$
for some $\te \in \R \setminus \Q.$
\end{prp}

\begin{proof}
Let $\Diff (M \times S^1)$ be the set of all $\Ci$ diffeomorphisms of
$M \times S^1,$
and, following Section~5 of~\cite{FH}, give it the $\Ci$ topology.
For $\ld \in S^1,$ define $r_{\ld} \colon M \times S^1 \to M \times S^1$
by $r_{\ld} (x, \zt) = (x, \ld \zt).$
This defines a free smooth action of $S^1$ on $M \times S^1.$
Now let $D \subset \Diff (M \times S^1)$ be the closure of the set
\[
D_0 = \ts{ \left\{ g \circ r_{\ld} \circ g^{-1} :
   g \in \Diff (M \times S^1), \, \ld \in S^1 \right\} }.
\]
Fathi and Herman show that the subset $S$ of $D$ consisting of those
elements which are minimal is a dense $G_{\dt}$-set in $D$
(5.6 of~\cite{FH}),
and that the subset $T$ of $D$ consisting of those
elements which are uniquely ergodic is a dense $G_{\dt}$-set in $D$
(6.5 of~\cite{FH}).
All elements of $D$ are homotopic to the identity map,
because the $r_{\ld}$ are.
Since $D$ can be given a complete metric, it therefore suffices
to show that the set of $h \in D,$ for which there is some invariant
Borel probability measure $\mu$ with
$\rh_h^{\mu} ([u]) \in S^1 \setminus \exp (2 \pi i \Q),$
is a dense $G_{\dt}$-set in $D.$

For $\ld \in S^1,$ let
$E_{\ld} \subset D$ consist of those $h \in D$ for which there is an
invariant
Borel probability measure $\mu$ with $\rh_h^{\mu} ([u]) = \ld.$
Lemma~\ref{RotNoSC} implies that $E_{\ld}$ is closed.
We show that its complement is dense.

Let $\om \in S^1.$
Write $\om = \exp (2 \pi i \af)$ with $\af \in \R.$
Let $a$ be the constant function $a (x, \zt) = \af$ for all
$(x, \zt) \in M \times S^1.$
Then
\[
u \ts{ \left( \rsz{ r_{\om}^{-1} } (x, \zt) \right)^* } u (x, \zt)
   = \exp (2 \pi i a (x, \zt)).
\]
For any invariant measure $\mu,$ we can then compute
$\rh_{r_{\om}}^{\mu} ([u])$ as
\[
\rh_{r_{\om}}^{\mu} ( [u] )
  = \exp \left( 2 \pi i \int_X a \, d \mu \right)
  = \exp (2 \pi i \af) = \om.
\]

Now consider $\rh_h^{\nu} ([u])$
with $h = g \circ r_{\om} \circ g^{-1}$
for some $g \in \Diff (M \times S^1).$
We can write
\[
\left[ \ts{ \left( \rsz{ u \circ g^{-1} } \right) } \ts{ \left( \rsz{
        g \circ r_{\om}^{-1} } \circ g^{-1} (x, \zt) \right)}\right]^*
     \ts{ \left( \rsz{ u \circ g^{-1} } \right) } (x, \zt)
   = \ts{ \exp \left( 2 \pi i \left( \rsz{ a \circ g^{-1} } \right)
          (x, \zt) \right) } = \om,
\]
so that $\rh_h^{\nu} ( [u \circ g^{-1}] ) = \om$ for any
Borel probability measure $\nu$ which is invariant under $h.$
It follows from Lemma~\ref{hty} that
$[u \circ g^{-1}] = [u]$ or $[u \circ g^{-1}] = [u^{-1}].$
Therefore $\rh_h^{\nu} ( [u] ) \in \{ \om, \, \om^{-1} \}.$
In particular, if $\om \not\in \{ \ld, \, \ld^{-1} \},$
then $g \circ r_{\om} \circ g^{-1} \not\in E_{\ld}.$
Since $\{ r_{\om} \colon \om \in S^1 \setminus \{ \ld, \, \ld^{-1} \} \}$
is clearly dense in $\{ r_{\om} \colon \om \in S^1 \},$
it follows that $D \setminus E_{\ld}$ is dense in $D,$ as claimed.

Since $D$ has a complete metric, the set
$\bigcap_{\ld \in \Q} (D \setminus E_{\ld} )$
is a dense $G_{\dt}$-set in~$D.$
Moreover, its intersection
\[
S \cap T \cap \bigcap_{\ld \in \Q} (D \setminus E_{\ld} )
\]
with the dense $G_{\dt}$-sets $S$ and $T$ from the beginning of
the proof is again a dense $G_{\dt}$-set, and in particular
not empty.
But any element of this set is a minimal diffeomorphism of
$M \times S^1$ which satisfies the conclusion of the proposition.
\end{proof}

It seems reasonable to expect that, for any compact manifold $M$
and any $\te \not\in \Q,$
there is a uniquely ergodic \md\  such that the rotation number
of $[u]$ with respect to its unique invariant measure is $\te.$
Proving this, however, requires more work.

\begin{exa}\label{Ex3}
Let $M_1 = S^2 \times S^1,$ and let $u \in U (C (M_1))$
be given by $u (x, \zt) = \zt.$
Use Proposition~\ref{Exist}
to choose $\te \in [0, 1] \setminus \Q$ and a 
uniquely ergodic minimal diffeomorphism $h_1 \colon M_1 \to M_1,$ with
unique invariant  Borel probability measure $\mu,$
which is homotopic to the identity map and such
that $\rh_{h_1}^{\mu} ([u]) = \exp (2 \pi i \te) \in S^1.$
Let $M_2 = (S^1)^3,$ and define
$h_2 \colon M_2 \to M_2$ by
\[
h_2 (\zt_1, \zt_2, \zt_3)
 = \left( \exp (2 \pi i \te) \zt_1,
             \, \zt_1 \zt_2, \, \zt_2 \zt_3 \right)
\]
for $(\zt_1, \zt_2, \zt_3) \in (S^1)^3.$
We show that
$C^* (\Z, M_1, h_1) \cong C^* (\Z, M_2, h_2).$
Note that $h_1$ can't be topologically orbit equivalent to $h_2,$
because $M_1$ is not homeomorphic to $M_2.$

Lemma~\ref{FT3D} implies that $h_2$ is minimal and uniquely ergodic,
and computes the Elliott invariant of $C^* (\Z, M_2, h_2).$
We calculate the Elliott invariant of $C^* (\Z, M_1, h_1).$

The K\"{u}nneth Theorem (Theorem~4.1 of~\cite{Sc1}; also see
Corollary 2.7.15 of~\cite{At} for the commutative case,
which suffices here)
shows that
\[
K^* (S^2 \times S^1) \cong K^* ( S^2) \otimes K^* ( S^1).
\]
We identify the two sides of this isomorphism.
Let $1 \in C (S^2)$ be the identity, and let $\bt \in K^0 (S^2)$
be the Bott element, which is of the form $[p] - [q]$ for rank one
\pj s $p, \, q \in M_2 (C (S^2)).$
Let $1 \in C (S^1)$ be the identity, and let $z \in C (S^1)$
be the canonical unitary $z (\zt) = \zt.$
Then $K^0 (S^2 \times S^1)$
is the free abelian group on generators
\[
\et_1 = [1] \otimes [1] \andeqn \et_2 = \bt \otimes [1],
\]
and $K^1 (S^2 \times S^1)$
is the free abelian group on generators
\[
\gm_1 = [1] \otimes [z] = [u] \andeqn \gm_2 = \bt \otimes [z].
\]
Since $h_1$ is homotopic to the identity map (by construction),
$h_1^* = \id,$ and the exact sequence of Theorem~\ref{PVSeq}
breaks apart into the two exact sequences
\[
0 \longrightarrow \Z \et_1 + \Z \et_2
  \longrightarrow K_0 ( C^* (\Z, M_1, h_1))
  \stackrel{\partial}{\longrightarrow} \Z \gm_1 + \Z \gm_2
  \longrightarrow 0
\]
and
\[
0 \longrightarrow \Z^2
  \longrightarrow K_1 ( C^* (\Z, M_1, h_1))
  \longrightarrow \Z^2
  \longrightarrow 0.
\]
Both split because $\Z^2$ is free.
Therefore
\[
K_0 ( C^* (\Z, M_1, h_1)) \cong \Z^4 \cong K_0 ( C^* (\Z, M_2, h_2))
\]
and
\[
K_1 ( C^* (\Z, M_1, h_1)) \cong \Z^4 \cong K_1 ( C^* (\Z, M_2, h_2)).
\]
We are done with $K_1 ( C^* (\Z, M_1, h_1) ),$
but it remains to determine
the action of the tracial state and the order on $K_0 ( C^* (\Z, M_1, h_1) ).$

Identify $\et_1$ and $\et_2$
with their images in $K_0 ( C^* (\Z, M_1, h_1) ).$
Let $\ta$ be the unique tracial state on $C^* (\Z, M_1, h_1),$ which comes
from the unique ergodic measure $\mu$ on $M_1.$
Clearly $\ta_* ( \et_1 ) = 1.$
We have $\ta_* ( \et_2 ) = 0$ because
$\bt = [p \otimes 1] - [q \otimes 1]$ and $p$ and $q$ have the same
rank at each point.

We have $K^1 (M_1)^{h_1} = K^1 (M_1)$ by the above, and
our next step is to calculate $\rh_h^{\mu} ( \gm_1 )$ and
$\rh_h^{\mu} ( \gm_2 )$ as in Definition~\ref{RDfn}.
We have $\rh_h^{\mu} ( \gm_1 ) = \exp (2 \pi i \te)$ by construction.
For $\gm_2,$ write $\bt = [p] - [q]$ as above.
Working in a suitable matrix algebra, we see that $\gm_2$ is represented
by the unitary
\[
[ (1 - p) \otimes 1 + p \otimes z ]
      [(1 - q) \otimes 1 + q \otimes z^{-1} ].
\]
Its determinant is the constant function $1,$
so $\rh_h^{\mu} ( \gm_2 ) = 1.$

Choose $\nu_1^{(0)}, \, \nu_2^{(0)} \in K_0 ( C^* (\Z, M_1, h_1) )$
such that
$\ts{ \partial \left( \rsz{ \nu_1^{(0)} } \right) } = \gm_1$
and
$\ts{ \partial \left( \rsz{ \nu_2^{(0)} } \right) } = \gm_2.$
Theorem~\ref{RThm} therefore implies that
\[
\ts{ \ta_* \left( \rsz{ \nu_1^{(0)} } \right) } \in \te + \Z
\andeqn
\ts{ \ta_* \left( \rsz{ \nu_2^{(0)} } \right) } \in \Z.
\]
Taking $\nu_1 = \nu_1^{(0)} - k \overline{\et}_1$ and
$\nu_2 = \nu_2^{(0)} - l \overline{\et}_1$ for suitable $k$ and $l,$
we get $\ta_* (\nu_1) = \te$ and $\ta_* (\nu_2) = 0.$

We can now identify $K_0 ( C^* (\Z, M_1, h_1) )$ as
\[
\Z \et_1 \oplus \Z \et_2 \oplus \Z \nu_1 \oplus \Z \nu_2,
\]
with
\[
\ta_* ( \et_1 ) = 1, \,\,\,\,\,\,
\ta_* ( \nu_1 ) = \te, \andeqn
\ta_* ( \et_2 ) = \ta_* ( \nu_4 ) = 0.
\]
Comparing this with the result of Lemma~\ref{FT3D} for the case
$m = n = 1,$ we see that there is an isomorphism
$f \colon K_0 ( C^* (\Z, M_1, h_1) ) \to K_0 ( C^* (\Z, M_2, h_2) )$
which preserves the tracial states and the class of the identity.
We already have the isomorphism on $K_1,$ and the range of the
tracial state on $K_0$ is dense, so Theorem~\ref{Isom} implies that
$C^* (\Z, M_1, h_1) \cong C^* (\Z, M_2, h_2).$

{}From Proposition~\ref{NonconjCartan},
we see that the algebra has a diagonal (Cartan subalgebra),
isomorphic to $C (M_1),$ and another one isomorphic to $C (M_2).$
\end{exa}

\begin{lem}
The \md s $h_1$ and $h_2$ of Example~\ref{Ex3} are not
flip flow equivalent.
\end{lem}

\begin{proof}
It follows from Theorem~2 of~\cite{Sz} that if $h_1$ and $h_2$
are \hme s of connected \cms s $X_1$ and $X_2,$
and if $h_1$ and $h_2$ are flip flow equivalent,
then the universal covers of $X_1$ and $X_2$ are homeomorphic.
Clearly, however, the universal covers of $S^2 \times S^1$ and
$(S^1)^3$ are not homeomorphic.
\end{proof}

\section{Fourth example: Minimal diffeomorphisms on manifolds of different
    dimensions}\label{Sec:DiffDims}

We find uniquely ergodic \md s  $h_n \colon S^n \times S^1 \to S^n \times S^1,$
for $n \geq 3$ odd, such that the \cp\  \ca s
$C^* (\Z, \, S^n \times S^1, \, h_n)$ are all isomorphic.
No two of these \diff s can be topologically conjugate, or even
flip flow equivalent, since they act on \mf s of different dimensions.
In terms of the three additional ways to think of the results:
\begin{itemize}
\item
The dimension of the space on which the \hme\  acts
is not an invariant of the \tgca.
\item
The \ca\  in the example
has infinitely many diagonals (Cartan subalgebras)
which are far from being isomorphic to each other,
having maximal ideal spaces of different dimensions.
\item
We don't know if the \md s in this example
can be chosen to be tempered in the sense of
Definition~3.1 of~\cite{Ph},
but it seems reasonable to hope that they can be.
If so, this example seems to be a good candidate for one in
which the smooth crossed products are not isomorphic.
\end{itemize}

The basis of the construction is the existence of
uniquely ergodic \md s on odd spheres and the following result,
obtained by combining several results from~\cite{Pr}.

\begin{lem}\label{RProd}
Let $X$ be a connected \cms,
and let $h \colon X \to X$ be a uniquely
ergodic \mh.
Then there is a dense $G_{\dt}$-set $T \subset S^1$ such that, for every
$\ld \in T,$ the \hme\  of $X \times S^1$ given by
$(x, \zt) \mapsto (h (x), \ld \zt)$ is minimal and uniquely ergodic.
\end{lem}

\begin{proof}
For $\ld \in S^1,$ define $r_{\ld} \colon X \times S^1 \to X \times S^1$
by $r_{\ld} (x, \zt) = (x, \ld \zt).$
This is a \ct\  action of $S^1$ on $X \times S^1.$
We claim that
the pair $( X \times S^1, \, h \times \id )$ is a simple free
extension of $(X, h)$ in the sense of~\cite{Pr} (see Equation~(0.4)).
Indeed, for $n \in \widehat{S^1} \cong \Z,$
the required \cfn\  $f_n \colon X \times S^1 \to S^1$ can be taken to be
simply $f_n (x, \zt) = \zt^n.$

Since $X \times S^1$ is connected, Theorem~2 of~\cite{Pr} implies that
\[
\{ \ld \in S^1 :
   r_{\ld} \circ (h \times \id) \,\, {\text{is minimal}} \}
\]
contains a dense $G_{\dt}$-set in $S^1.$
Since $h$ is uniquely ergodic, we use the version of
Theorem~4 of~\cite{Pr} indicated in Remark~3 (following Theorem~5) there
to show that
\[
\{ \ld \in S^1 :
   r_{\ld} \circ (h \times \id) \,\, {\text{is uniquely ergodic}} \}
\]
contains a dense $G_{\dt}$-set in $S^1.$
The set $T$ of the lemma is obtained by intersecting these two sets.
\end{proof}

\begin{exa}
Use Theorem~3 (in Section~3.8) of~\cite{FH} to find,
for each odd $n \geq 3,$
a uniquely ergodic \md\  $h_n^{(0)} \colon S^n \to S^n.$
We note that $h_n^{(0)}$ must have degree $1,$ since the
Lefschetz fixed point theorem (Theorem 4.7.7 of~\cite{Sp})
implies that an orientation
reversing \diff\  of an odd sphere must have a fixed point;
see 4.7.9 of~\cite{Sp}.
Therefore $h_n^{(0)}$ is homotopic to the identity map.

Let $T_n \subset S^1$ be the dense $G_{\dt}$-set associated with
$h_n^{(0)}$ from Lemma~\ref{RProd}.
Then $T = \bigcap_{k = 1}^{\infty} T_{2 k + 1}$ is still
a dense $G_{\dt}$-set.
Choose $\te \in [0, 1]$ such that $\exp (2 \pi i \te) \in T.$
Define $h_n \colon S^n \times S^1 \to S^n \times S^1$ by
$h_n (x, \zt)
 = \left( \rsz{ h_n^{(0)} (x) }, \, \exp ( 2 \pi i \te) \zt \right).$
Then each $h_n$ is minimal and uniquely ergodic, by the choice of $T.$
In particular, $\te \not\in \Q.$
Also, each $h_n$ is homotopic to the identity map.
We prove that the \ca s
$C^* (\Z, \, S^n \times S^1, \, h_n)$ are all isomorphic.
No two of these \diff s can be topologically conjugate, or even
flip flow equivalent, since they act on \mf s of different dimensions.

As usual, we start by using the
K\"{u}nneth Theorem (Theorem~4.1 of~\cite{Sc1} or
Corollary 2.7.15 of~\cite{At}), to get
\[
K^* (S^n \times S^1) \cong K^* ( S^n) \otimes K^* ( S^1),
\]
and we identify the two sides of this isomorphism.
We write
\[
K^0 (S^n) = \Z \cdot [1], \,\,\,\,\,\,
K^1 (S^n) = \Z \cdot \gm, \,\,\,\,\,\,
K^0 (S^1) = \Z \cdot [1], \andeqn
K^1 (S^1) = \Z \cdot [z],
\]
with $\gm$ equal to the class of a suitable unitary in some
matrix algebra over $C (S^n)$ and with $z \in C (S^1)$
given by $z (\zt) = \zt.$
Then, suppressing the dependence on $n$ in the notation,
 $K^0 (S^n \times S^1)$
is the free abelian group on generators
\[
\et_1 = [1] \otimes [1] \andeqn \et_2 = \gm \otimes [z],
\]
and $K^1 (S^n \times S^1)$
is the free abelian group on generators
\[
\gm_1 = [1] \otimes [z] \andeqn \gm_2 = \gm \otimes [1].
\]
The exact sequence of Theorem~\ref{PVSeq}
breaks apart into the two exact sequences
\[
0 \longrightarrow \Z \et_1 + \Z \et_2
  \longrightarrow K_0 ( C^* (\Z, \, S^n \times S^1, \, h_n) )
  \stackrel{\partial}{\longrightarrow} \Z \gm_1 + \Z \gm_2
  \longrightarrow 0
\]
and
\[
0 \longrightarrow \Z^2
  \longrightarrow K_1 ( C^* (\Z, \, S^n \times S^1, \, h_n) )
  \longrightarrow \Z^2
  \longrightarrow 0.
\]
Both split because $\Z^2$ is free, so (with obvious identifications,
and with $\ts{ \partial \left( \rsz{ \nu_j^{(0)} } \right) } = \gm_j$)
\[
K_0 ( C^* (\Z, \, S^n \times S^1, \, h_n) ) =
\Z \et_1 + \Z \et_2 + \Z \nu_1^{(0)} + \Z \nu_1^{(0)}
\]
and
\[
K_1 ( C^* (\Z, \, S^n \times S^1, \, h_n) ) \cong \Z^4.
\]

Let $\ta$ be the unique tracial state on
$C^* (\Z, \, S^n \times S^1, \, h_n),$
which comes from the product $\mu \times \ld$ of the
unique ergodic measure $\mu$ for $h_n^{(0)}$ on $S^n$ and
normalized Haar measure $\ld$ on $S^1.$
(This product is clearly $h_n$-invariant, and by construction there is
only one invariant probability measure.)
Clearly $\ta_* ( \et_1 ) = 1.$
Since $\et_2$ vanishes under point evaluations (by naturality
in the K\"{u}nneth Theorem), we have $\ta_* ( \et_2 ) = 0.$
We have $K^1 (M_1)^{h_1} = K^1 (M_1)$ by the above, and
our next step is to calculate $\rh_{h_n}^{\mu \times \ld} ( \gm_1 )$ and
$\rh_{h_n}^{\mu \times \ld} ( \gm_2 )$ as in Definition~\ref{RDfn}.
It is easy to check
(using $n \geq 3$ and  Corollary VI.12(i) of~\cite{Ex}
for the second step) that
\[
\rh_{h_n}^{\mu \times \ld} ( \gm \otimes [1])
   = \rh_{h_n^{(0)}}^{\mu} ( \gm ) = 1,
\]
and (with $r$ being rotation by $\exp (2 \pi i \te)$) that
\[
\rh_{h_n}^{\mu \times \ld} ( [1] \otimes [z]) = \rh_r^{\ld} ( [z] )
   = \exp (2 \pi i \te).
\]
As usual, we define $\nu_1$ and $\nu_2$ by subtracting suitable
multiples of $[1]$ from $\nu_1^{(0)}$ and $\nu_2^{(0)}.$
This gives
\[
K_0 ( C^* (\Z, \, S^n \times S^1, \, h_n) ) =
\Z \et_1 + \Z \et_2 + \Z \nu_1 + \Z \nu_1
\]
with
\[
[1] = \et_1, \,\,\,\,\,\,
\ta_* ( \et_1 ) = 1, \,\,\,\,\,\,
\ta_* ( \nu_1 ) = \te, \andeqn
\ta_* ( \et_2 ) = \ta_* ( \nu_4 ) = 0.
\]
The generators $\et_1, \, \et_2, \, \nu_1, \, \nu_1$ depend on $n,$
but it is clear from this formula that for any two values of $n$
there is an isomorphism between these groups which preserves
$\ta_*$ and~$[1].$
Since the range of the tracial state on $K_0$ is dense,
Theorem~\ref{Isom} implies that the \ca s are pairwise isomorphic.

{}From Proposition~\ref{NonconjCartan},
we see that an algebra in the common isomorphism class
has, for every odd $n \geq 3,$ a diagonal (Cartan subalgebra)
isomorphic to $C (S^n \times S^1).$
\end{exa}

\end{document}